\documentclass[conference]{IEEEtran}
\IEEEoverridecommandlockouts
\usepackage{cite}
\usepackage{amsmath,amssymb,amsfonts}
\usepackage{algorithmic}
\usepackage{graphicx}
\usepackage{textcomp}
\usepackage{xcolor}
\def\BibTeX{{\rm B\kern-.05em{\sc i\kern-.025em b}\kern-.08em
    T\kern-.1667em\lower.7ex\hbox{E}\kern-.125emX}}

\usepackage{amsthm}
\usepackage{float}
\usepackage{commath}
\usepackage{mathtools}
\usepackage{nicefrac}
\usepackage{hyperref}
\usepackage{algorithm2e}

\newcommand{\cola}{\textsc{CoLa}}%

\newcommand{\grad}{\nabla}

\newcommand{\Rr}{\mathbb{R}}
\newcommand{\Rrv}[1]{\Rr^{#1}}

\newcommand{\inner}[2]{#1 \cdot #2}

\newcommand{\vect}[1]{\mathbf{#1}}
\newcommand{\bb}{\vect b}

\newcommand{\vv}{\vect v}

\newcommand{\xx}{\vect x}
\newcommand{\yy}{\vect y}
\newcommand{\zz}{\vect z}

\newcommand{\Aa}{\vect A}
\newcommand{\Ak}{\vect{A}_{[k]}}
\newcommand{\xk}{\vect{x}_{[k]}}
\newcommand{\dxk}{\Delta\vect{x}_{[k]}}

\newcommand{\vk}{\vect{v}_{k}}

\newcommand{\dvk}{\Delta\vect{v}_{k}}
\newcommand{\Gk}{\Gamma_{k} (\dxk; \vk, \xk)}
\newcommand{\gk}{g_{[k]}}

\newcommand{\Oa}{\mathcal{O}_A}

\newcommand{\ootau}{$(\nicefrac{1}{\tau})-$}

\begin{document}

\title{Distribution System Voltage Prediction from Smart Inverters using Decentralized Regression\\
\thanks{This work is supported by the U.S. Department of Energy’s Office of Energy Efficiency and Renewable Energy (EERE) under Solar Energy Technologies Office (SETO) Agreement Number 34229.  Prepared by LLNL under Contract DE-AC52-07NA27344.  LLNL-CONF-816452.}
}

\author{\IEEEauthorblockN{Zachary R. Atkins}
\IEEEauthorblockA{\textit{Department of Mathematics} \\
\textit{University of Kansas}\\
Lawrence, USA \\
zatkins@ku.edu}
\and
\IEEEauthorblockN{Christopher J. Vogl}
\IEEEauthorblockA{\textit{Center for Applied Scientific Computing} \\
\textit{Lawrence Livermore National Laboratory}\\
Livermore, USA \\
vogl2@llnl.gov}
\and
\IEEEauthorblockN{Achintya Madduri}
\IEEEauthorblockA{\textit{Computational Engineering Division} \\
 \textit{Lawrence Livermore National Laboratory}\\
 Livermore, USA \\
 madduri2@llnl.gov}
\and
\IEEEauthorblockN{Nan Duan}
\IEEEauthorblockA{\textit{Computational Engineering Division} \\
\textit{Lawrence Livermore National Laboratory}\\
Livermore, USA \\
duan4@llnl.gov}
\and
\IEEEauthorblockN{Agnieszka K. Miedlar}
\IEEEauthorblockA{\textit{Department of Mathematics} \\
\textit{University of Kansas}\\
Lawrence, USA\\
amiedlar@ku.edu}
\and
\IEEEauthorblockN{Daniel Merl}
\IEEEauthorblockA{\textit{Center for Applied Scientific Computing} \\
\textit{Lawrence Livermore National Laboratory}\\
Livermore, USA\\
dmerl@llnl.gov}
}

\maketitle

\begin{abstract}
  As photovoltaic (PV) penetration continues to rise and smart inverter functionality continues to expand, smart inverters and other distributed energy resources (DERs) will play increasingly important roles in distribution system power management and security.  In this paper, it is demonstrated that a constellation of
  smart inverters in a simulated distribution circuit can enable precise voltage
  predictions using an asynchronous and decentralized prediction algorithm.  Using simulated data and a
  constellation of $15$ inverters in a ring communication topology, the
  \cola{} algorithm is shown to accomplish the learning task required for voltage
  magnitude prediction with far less communication overhead than fully connected
  P2P learning protocols.  Additionally, a dynamic stopping criterion is proposed that does not require a regularizer like the original \cola{} stopping criterion.
\end{abstract}


\section{Introduction}
Power distribution networks are increasingly dynamic due to the growing
prevalence of distributed energy resources (DERs), such as rooftop solar photovoltaic (PV) panels.
The old distribution network model as a predictable set of loads is no longer accurate.
In particular, it is not uncommon for distribution networks with high PV-penetration to be net producers of power during certain parts of the year.
As a result, distribution networks with high PV-penetration require more complex monitoring and
control objectives compared to traditional networks \cite{wang_inverter_2020}.  Meeting
the demands of complex monitoring requires large-scale deployment of smart
sensors and efficient processing of large volumes of data
\cite{sanduleac_meter_2020}.  This has motivated recent research into utilizing
``on-device'' computation to avoid collecting data from all smart devices in one
location for processing \cite{duan_inference_2019}.  The industry standard today
is to collect sensor data from the smart devices in a central control facility,
where it is analyzed for the purposes of maintaining a healthy and safe grid.
Collecting data in one location, however, can incur large communication cost and
can require large data storage capability \cite{meng_distributed_2018}.
Additionally, there can be cybersecurity concerns (e.g. data breaches)
\cite{soyoye_risk_2019} and a host of privacy concerns ranging from intellectual
property loss \cite{sakurama_decentralized_2017} to personal data loss
\cite{dep2sa}.

Without addressing these issues, power utilities are unable to safely accept
more solar generation in the grid.  If the computational resources of smart grid
devices, such as smart meters and smart solar inverters can be exploited, the
stress on the central control facility can be alleviated by migrating some
DER analyses and control decisions to either distributed energy resource
management system (DERMS) controllers or the smart inverters themselves
\cite{saxena_dds_2018}.  Also, in the event that a cyber attack targeting smart
devices is launched \cite{duan_cyber_2020}, computational resources of smart
devices can facilitate local command verification to improve the cybersecurity
posture of the power grid \cite{zografopoulos_derauth_2020}.  Therefore, the
development of ``on-device'' algorithms is essential for power utilities to
harness the full potential of DER power generation.

The focus of this work is on the optimization/learning/training of a predictive
model across edge devices, such as smart inverters, for the purpose of assessing the health of the grid (e.g. over-voltage
issues).  There is a vast array of techniques from the optimization and
machine learning communities for accomplishing distributed linear regression.
That said, the computing and data transfer limitations of a collection of edge
devices requires a distributed algorithm that also avoids (i) a central
controller/aggregator (decentralized) and (ii) global synchronization of
local data (asynchronous). The \textit{\textbf{CO}mmunication-Efficient
Decentralized \textbf{L}inear Le\textbf{A}rning} (\cola{}) algorithm proposed in
\cite{he_cola_2018} is specifically designed for
``column-partitioned data'', in which each edge device represents a sort of data
silo from which data egress is to be minimized.

This work investigates the applicability of the \cola{} algorithm to linear regression training on smart inverter data.  This work
also suggests a stopping criterion for the \cola{} algorithm based on the
magnitude of the local iterate updates. Such a stopping criterion enables a
broader class of global objective functions than the ones introduced in \cite{he_cola_2018}, e.g., those that do not require a regularizer function.
The remainder of this paper is organized as follows: Section \ref{sec:learning} describes the \cola{} algorithm and proposes a dynamic stopping
criterion; Section \ref{sec:regression} formulates the inverter data
regression problem; Section \ref{sec:results} presents the numerical results of
applying the \cola{} algorithm to inverter voltage magnitude data; finally
Section \ref{sec:conclusion} draws conclusions.

\section{Decentralized Learning Method \label{sec:learning}}
Consider a collection of $n$ smart inverters collecting local voltage magnitude at regular intervals (e.g. every $15$ minutes for a day), resulting in $m$ voltage magnitude values for each inverter.
Denote as $\Aa \in \mathbb{R}^{m \times n}$ the data matrix with each column representing voltage magnitude values from a single inverter.
Consider a smart meter also collecting local voltage magnitude at the same times as the smart inverters.
Denote as $\bb \in \mathbb{R}^m$ the data vector containing the smart meter voltage magnitudes values.
A regression model, such as linear least-squares, can be trained on $\Aa$ and $\bb$ and used to predict future smart meter voltage magnitude values from the local smart inverter data.
Let $f: \mathbb{R}^m \to \mathbb{R}$ be the objective/loss function to be minimized to train the regression model.
With $g: \mathbb{R}^n \to \mathbb{R}$ as a convex regularization function, the training of the linear regression is modeled as an optimization problem:
\begin{gather}
  \min_{\xx \in \mathbb{R}^n} \, f(\Aa \xx) + g(\xx).
  \label{eqn:minimization}
\end{gather}

\subsection{The \cola{} Algorithm}

  The \cola{} algorithm \cite{he_cola_2018} will be used to solve \eqref{eqn:minimization} in a decentralized fashion across $K$ nodes.
  The algorithm requires regularizer $g$ to be additively separable, $g(\xx) = \sum_{i=1}^n g_i(x_i)$, and objective function $f$ to be of the form $f(\vv)$, where $\vv = \Aa \xx$  is convex, and \ootau smooth, i.e.,
  \begin{gather*}
    f(\yy) \leq f(\zz) + \inner{\grad f(\zz)}{(\mathbf{y}-\mathbf{z})} + \frac{1}{2\tau}\norm{\mathbf{y}-\mathbf{z}}_2^2 \,\, \forall \, \yy, \zz \in \Rrv n.
  \end{gather*}
  The algorithm, shown as Algorithm \ref{alg:cola}, decomposes the global minimization problem \eqref{eqn:minimization} into a set of local subproblems, which are then solved independently of one another.
  The local subproblem for node $k$ is defined as
  \begin{gather}
    \begin{gathered}
      \Gk = \frac{1}{K} f(\vk) + \inner{\grad f(\vk)}{\Ak\dxk} + \\
      \frac{K}{2\tau}\norm{\Ak\dxk}_2^2 + \sum_{i \in \mathcal{P}_k} g_i(x_i + \Delta x_{[k]i}),
    \end{gathered}
    \label{eqn:subproblem}
  \end{gather}
  where $\mathcal{P}_k$ is a set of indices assigned to node $k$ that corresponds to the columns of $\Aa$ (denoted $\Ak$) and rows of $\xx$ (denoted $\xk$).
  The results from the local subproblems are mixed according to a doubly-stochastic mixing matrix $\mathcal{W} \in \mathbb{R}^{K \times K}$ in an asynchronous manner, making the \cola{} algorithm fully decentralized.
  For additional details on the relationship between solving the global minimization problem \eqref{eqn:minimization} and solving the local subproblems \eqref{eqn:subproblem}, see \cite{he_cola_2018}.
  \begin{algorithm}
    \caption{\cola{}: \textbf{Co}mmunication-Efficient Decentralized \textbf{L}inear Le\textbf{A}rning ($\gamma=1$).}\label{alg:cola}
    \SetAlgoLined\DontPrintSemicolon\LinesNumbered
    \SetKwInOut{Input}{Inputs}\SetKwInOut{Output}{Output}
    \SetKwFor{ParFor}{for}{in parallel do}{end}
    \SetKw{Update}{Update}%
    \SetKw{Init}{Initialize:}%
    \SetKw{Compute}{Compute}
    \Input{Data matrix $\Aa$ partitioned column-wise into $\{\Ak\}_{k=1}^K$, mixing matrix $\mathcal W$, and iteration number $N>0$.}
    \Output{Solution $\xx\in\Rrv n$ to problem \eqref{eqn:minimization}}
    \Init{$\xx^{(0)}\coloneqq \vect 0\in\Rrv n$; $\vk^{(0)}\coloneqq \vect 0\in\Rrv m\ \ \forall k=1\dots K$}\;
    \For{$i=0\dots N$}{
      \ParFor{$k=1\dots K$}{
        \Update{$\vk^{(i+\frac12)} = \sum_{\ell=1}^K \mathcal{W}_{kl}\vv_\ell^{(i)}$}\;
        \Compute{$\dxk\leftarrow$ solve a local subproblem $\Gk$}\;
        \Update{$\xk^{(i+1)} \coloneqq \xk^{(i)} + \dxk$}\;
        \Compute{$\dvk \coloneqq \Ak \dxk$ as a local estimate update}\;
        \Update{a local estimate $\vk^{(i+1)}=\vk^{(i+\frac12)} + K\dvk$}\;
      }
    }
  \end{algorithm}

\subsection{Stopping Criterion}
\label{sect:con:centering}

  As stated, Algorithm \ref{alg:cola} will stop after $N$ iterations regardless of whether the iterates have converged. In the case of centralized iterative algorithms, this static iteration count is usually replaced by a dynamic stopping criterion that may or may not depend on the current residual.
  Because the residual is a global quantity that usually requires synchronization of all nodes to be computed, a major challenge for decentralized iterative algorithms is determining a suitable stopping criterion that can be obtained in a decentralized fashion.
  The criterion proposed in \cite{he_cola_2018} utilizes a mapping between local data and a bound on the global duality gap.
  The local data in that mapping includes evaluating the complex conjugate of $g_i$, which excludes the use of $g(\xx) = 0$ and its corresponding complex conjugate $g^*(\xx^*) = \sup_{\xx^*} \, \inner{\xx^*}{\xx}$.
  As such, a decentralized stopping criterion is sought that indicates progress on minimizing $f(\Aa\xx)$ without requiring $g^*(\xx^*)$.

  Motivated by the distributed algorithms that use the magnitude of $\dxk$, the relationship between $\dxk$ and $f(\Aa\xx)$ is studied, i.e., the relationship between the backward and forward error.
  It is found that centering and normalizing the columns of data matrix,
  \begin{gather}
    \Aa_{ij} \rightarrow \frac{\Aa_{ij} - \mu_j}{\nu_j} \mbox{, \ where }
    \left \{
    \begin{aligned}
      \mu_j =& \frac{1}{m}\sum_{l=1}^m \Aa_{lj} \\
      \nu_j =& \left(\frac{1}{m} \sum_{l=1}^m \Aa_{lj}^2 \right)^{1/2}
    \end{aligned}
    \right .,
    \label{eqn:preprocess}
  \end{gather}
  have a significant impact on whether decreases in $\dxk$ correspond to decreases in $f(\Aa\xx)$.
  Because the data in the \cola{} algorithm is column-partitioned, \eqref{eqn:preprocess} can be applied in a fully decentralized fashion: no communication between nodes is necessary to preprocess the data.
  With the proposed preprocessing, the empirical results herein indicate that the magnitude of $\dxk$ is a viable candidate for a stopping criterion.

\section{Inverter Regression \label{sec:regression}}
\begin{figure}[!htbp]
  \centering
  \includegraphics[width=0.8\linewidth]{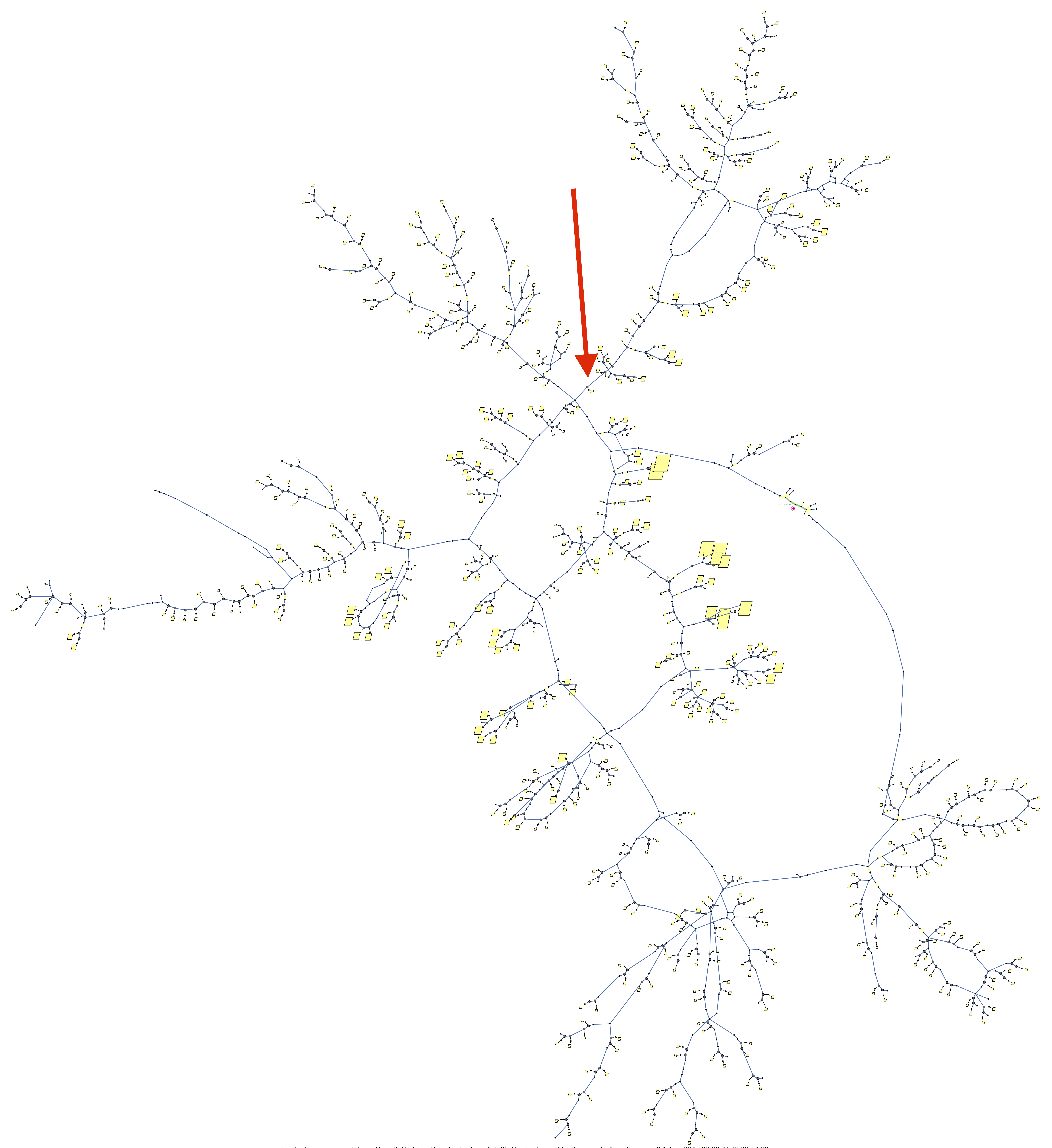}
  \caption[grid topology]{Network topology of the power distribution network used for the experiments.
  This distribution network has a peak demand of approximately $25$MW. Solar panels for a uniform $80\%$ grid penetration are indicated as yellow parallelograms with sizes represented relative to their peak capacity.
  The red arrow indicates the node where the voltage is being predicted by the regression-based algorithm.}
  \label{fig:grid_graph}
\end{figure}

 In this section, the utility of a regression-based voltage
prediction approach in simulated high-DER-penetration scenarios is demonstrated.  To obtain such
a network, the base model of the power distribution network shown in Figure
\ref{fig:grid_graph} is modified to include solar panels at $80\%$ penetration
levels. An associated solar panel was added to each load on the system such that
the peak power of the solar panels was equal to $80\%$ of the peak demand of the
associated load. The open source power grid simulator GridLAB-D
\cite{chassin_gridlabd_2008} is used to obtain the corresponding voltage data.
Given that only one distribution network is considered, the results herein are
presented as a case study in voltage prediction using a linear regression model.

\begin{figure}[!htbp]
  \centering
  \includegraphics[width=0.9\linewidth]{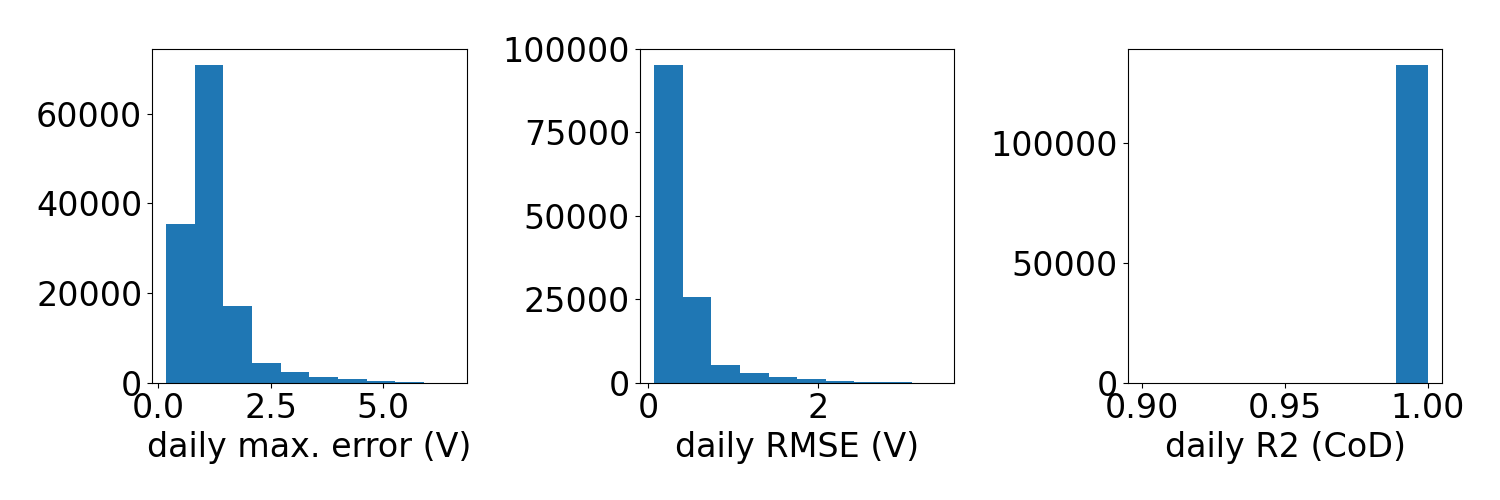}
  \caption[year regression results]{Linear least-square regression testing statistics, where each day is individually used to train a model that is testing on the remaining $364$ days (nominal voltage is 7200V).}
  \label{fig:yearly_regression}
\end{figure}

Because the node voltages in a high-PV-penetration distribution network can be
heavily impacted by the power generated by the distributed solar panels, the
impact of the particular training time periods on the ability of the algorithm
to predict node voltages accurately is first investigated.  The voltage
magnitude data obtained from $5$ inverters that are downstream from the node of
interest are used to train a linear least-squares regression model.  The
regression is trained on data obtained at $15$ minute intervals for a given day
of the year and then used to predict for all other days of the simulated year,
where each simulated day involves real stochastic variation in solar irradiation
levels, leading to variable PV generation.  The process is repeated so that each
calendar day is considered for training.  The results in Figure
\ref{fig:yearly_regression} show that linear least-squares regression model
trained on any day of the year can predict the voltage at the node of interest
for any other day with a high degree of accuracy.

With some confidence in the applicability of linear least-squares regression to
voltage prediction, the focus moves to using the \cola{} algorithm to train the
regression.  Voltage magnitude data is now obtained from $16$ solar inverters at
$15$ minute intervals for a given day of the yar.  A linear least-squares
regression model is used to predict values in $\bb$ from the data in $\Aa$:
\begin{gather*}
  f(\Aa \xx) = \frac{1}{2}[\inner{(\Aa \xx - \bb)}{(\Aa \xx - \bb)}],
\end{gather*}
where $\Aa$ is all but one column of the dataset, and $\bb$ is the final column.
It can be shown that $f$ is \ootau smooth, with $\tau = 1$, by noting
\begin{gather*}
  \frac{1}{2}\inner{(\yy - \zz)}{(\yy - \zz)} = \frac{1}{2} \inner{\yy}{\yy} - \frac{1}{2} \inner{\zz}{\zz} - \inner{\zz}{(\yy-\zz)} \\
  \Rightarrow \frac{1}{2}\norm{\yy}_2^2 = \frac{1}{2}\norm{\zz}_2^2 + \inner{\zz}{(\yy-\zz)} + \frac{1}{2} \norm{\yy - \zz}_2^2.
\end{gather*}
The local subproblem \eqref{eqn:subproblem} for the linear least-squares $f$  can be expressed as
\begin{gather*}
  \Gk = \frac{1}{2}\norm{\frac{\vk - \bb}{K} + \Ak\dxk}_2^2 \\
  + \gk(\xk + \dxk),
\end{gather*}
where each $\Ak$ is the column of $\Aa$ containing data from inverter $k$.

The implementation of the \cola{} algorithm provided in \cite{CoLA:repo} requires a non-zero regularizer $g$, which is a consequence of the choice of stopping criterion that uses the global duality gap in \cite{he_cola_2018}.
To satisfy this requirement, an elastic net regularizer is used:
\begin{gather*}
  g(\xx) = \lambda \left[ \frac{\eta}{2}(\inner{\xx}{\xx}) + (1-\eta) \sum_{j=1}^m |x_j| \right],
\end{gather*}
where $\eta\in(0,1]$ and $\lambda > 0$ are adjustable parameters.
As mentioned before, a ring topology network is used where each node $k$ exchanges $\vk$ values with only two other nodes.  In other words, $\mathcal{W}$ is symmetric with only three non-zero entries of $1/3$ in each row/column, with two of those entries being off-diagonal.
Another approach would be to have each of the $15$ ($n$) nodes broadcast the local data to all other nodes, so that each node has access to the full data matrix $\Aa$.
The latter approach involves the communication of $210$ ($n(n-1)$) column vectors of data.
With the ring topology, each \cola{} iteration involves the communication of $30$ ($2n$) column vectors.
Thus, if an acceptable regression is found in less than $7$ ($(n-1)/2$) iterations, the \cola{} algorithm results in less communication.
The specific implementation of the \cola{} algorithm used here is available in \cite{zenodo:plots}.

\section{Results \label{sec:results}}
\begin{figure}[!htbp]
  \centering
  \includegraphics[width=0.49\linewidth]{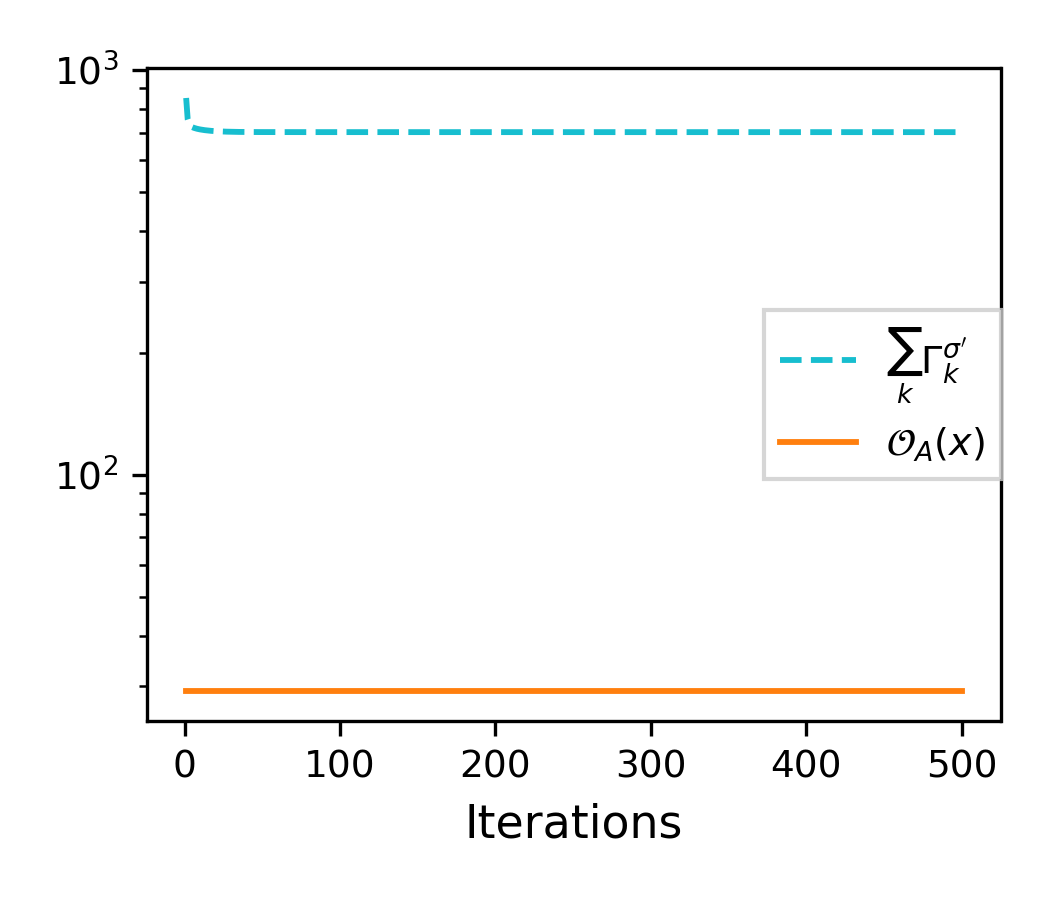}
  \includegraphics[width=0.49\linewidth]{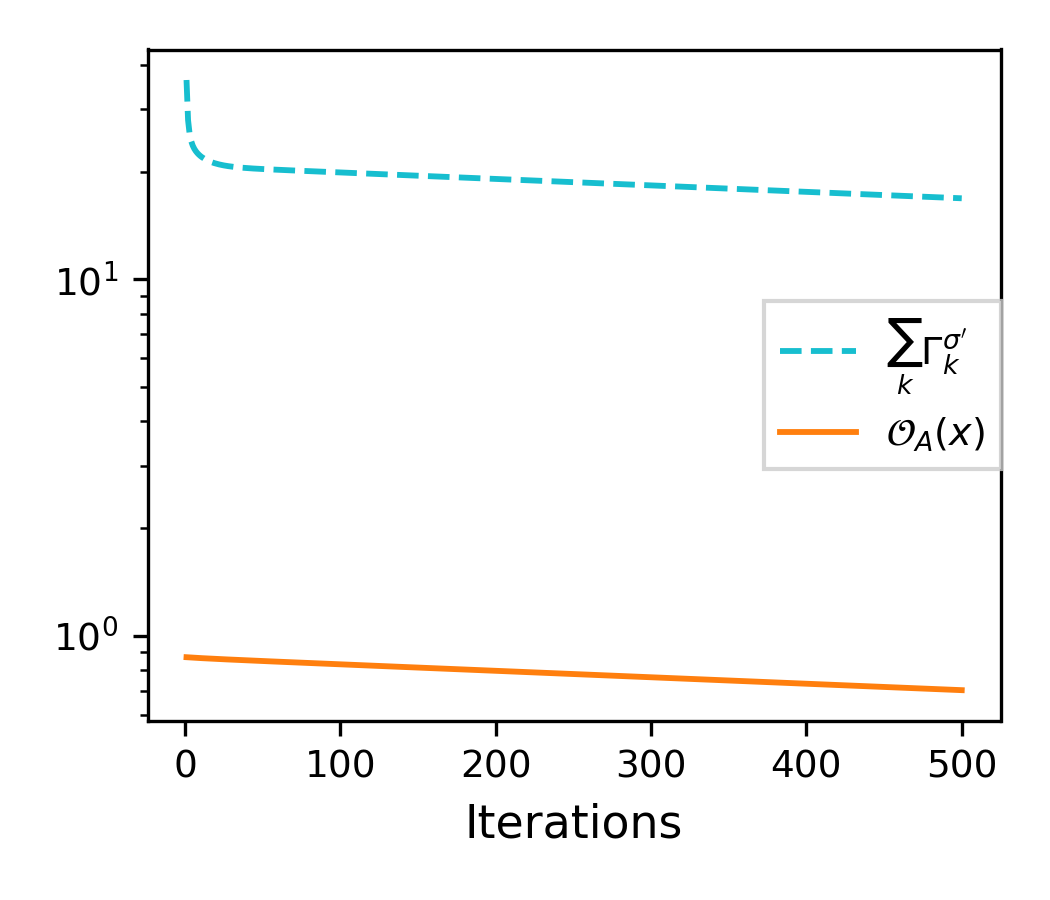}
  \caption[Local bounds on global objective for the inverter dataset]{Local subproblems and global objective function values for the inverter voltage dataset (left: original data, right: preprocessed data).}
  \label{fig:inv-thm1}
\end{figure}

To highlight the importance of data preprocessing, the relationship between the global minimization problem \eqref{eqn:minimization} and
solving the local subproblems \eqref{eqn:subproblem} is first observed. Figure \ref{fig:inv-thm1} demonstrates that the sum of the local subproblem objective functions $\Gk$ does indeed bound the value of the global objective function $\Oa(\xx) = f(\Aa \xx) + g(\xx)$. 
The decrease in the sum of $\Gk$ values in Figure \ref{fig:inv-thm1} is mirrored by the decrease in the values of $\Oa$, regardless of whether the inverter data is preprocessed or not.
That said, the preprocessing of the inverters voltage dataset appears to be necessary to avoid the stagnation of $\Oa$ values.

\begin{figure}[!htbp]
  \centering
  \includegraphics[width=0.9\linewidth]{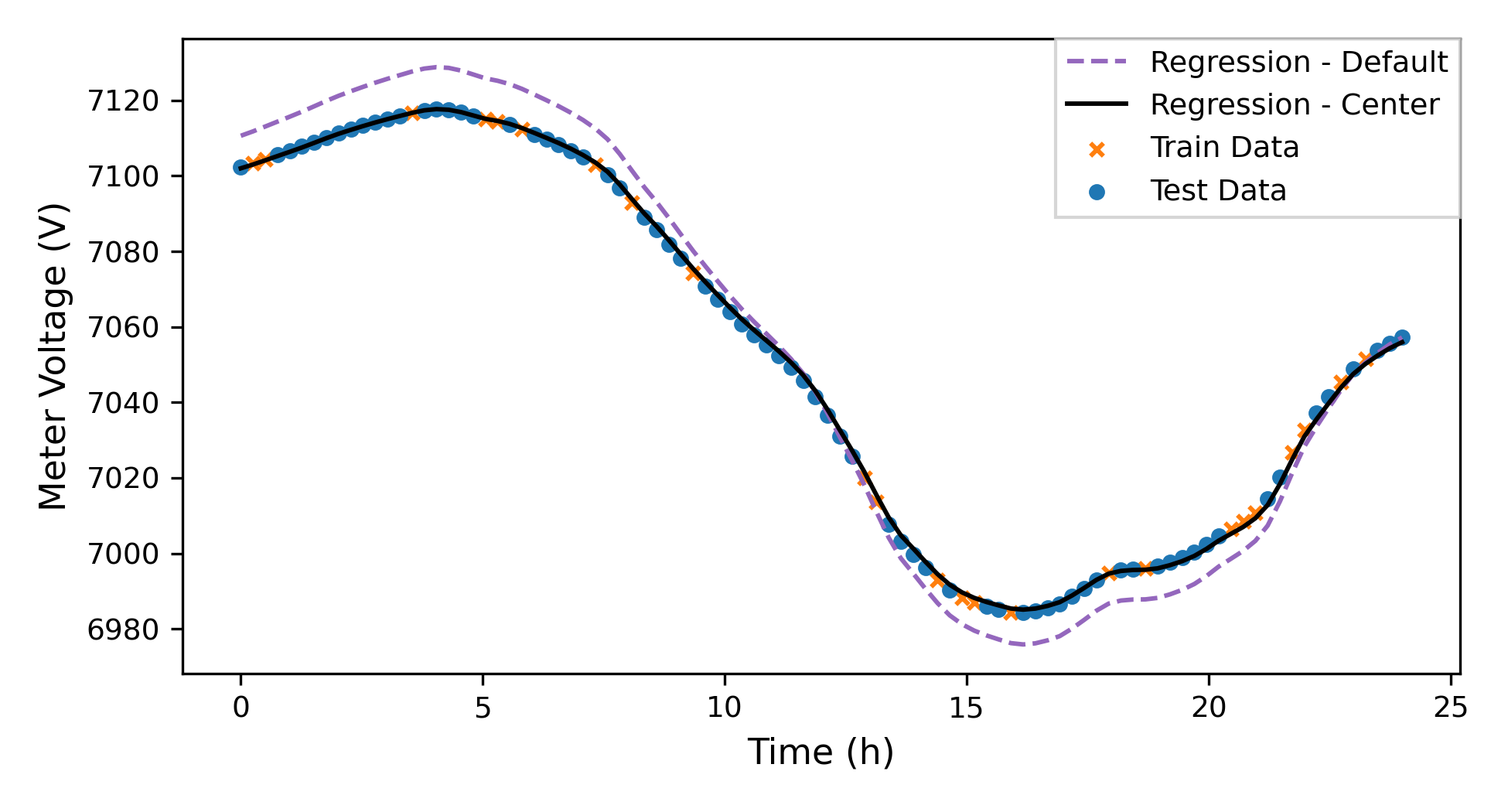}
  \caption[Regression comparison for the inverter dataset]{Regression comparison after $500$ iterations of the \cola{} algorithm applied to the elastic net problem with the inverter voltage data.}
  \label{fig:inv-regression}
\end{figure}

The stagnation in $\Oa$ values, when no preprocessing is performed, is reflected in the resulting regression quality obtained using the original data versus preprocessed data.
Figure \ref{fig:inv-regression} shows voltage magnitude profiles from both regressions after $500$ iterations, along with the corresponding train and test data.
The stagnation of $\Oa$ values in the case of the original data is reflected in a voltage magnitude profile that overshoots the peak voltage just before the $5^{th}$ hour (5 am) and undershoots the voltage dip around the $16^{th}$ hour (4 pm). On the other hand, the voltage magnitude profile that corresponds to the regression using preprocessed data very accurately captures both the corresponding train and test data.

\begin{figure}[!htbp]
    \centering
    \includegraphics[width=0.9\linewidth]{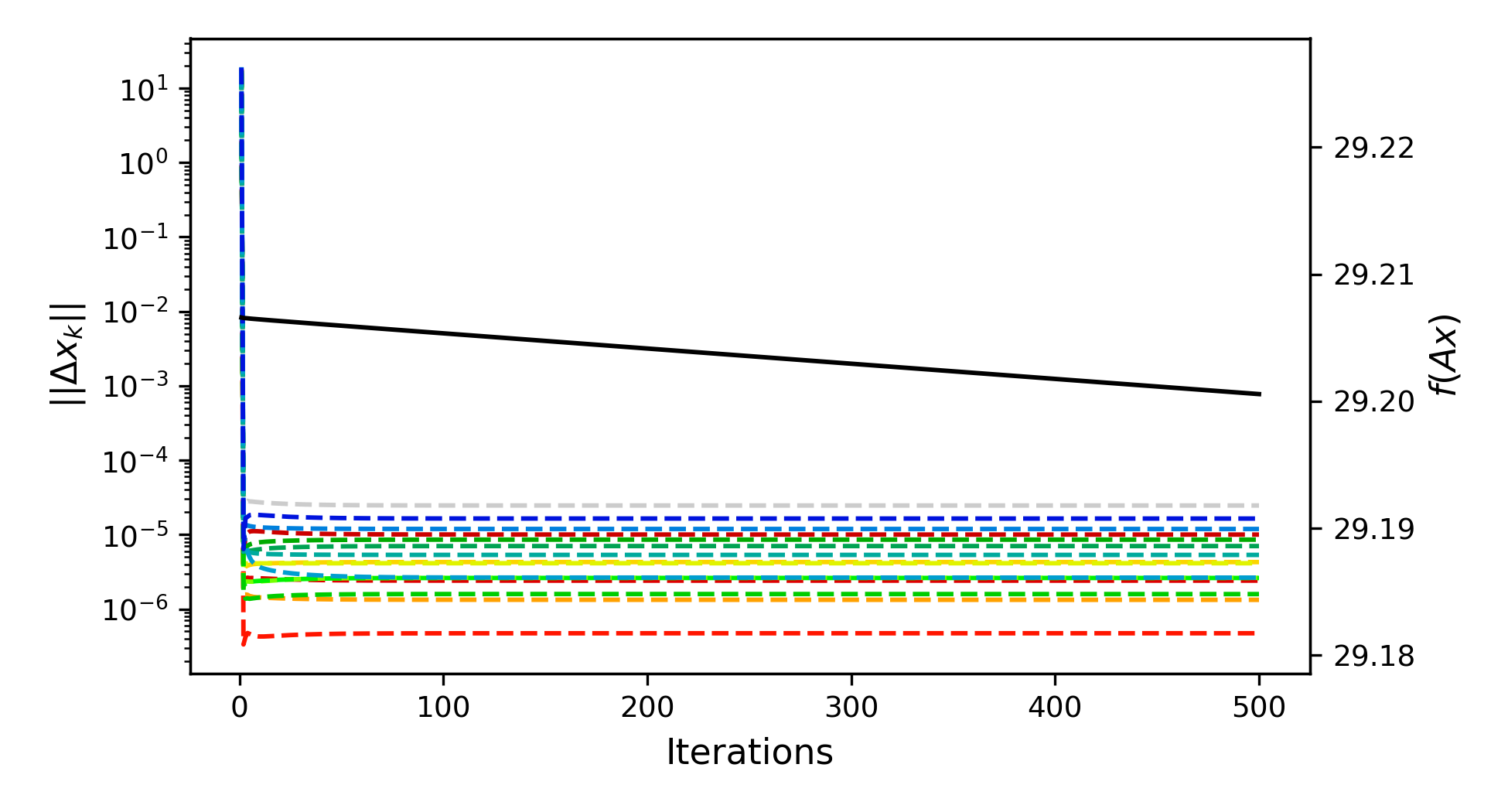}
    \\
    \includegraphics[width=0.9\linewidth]{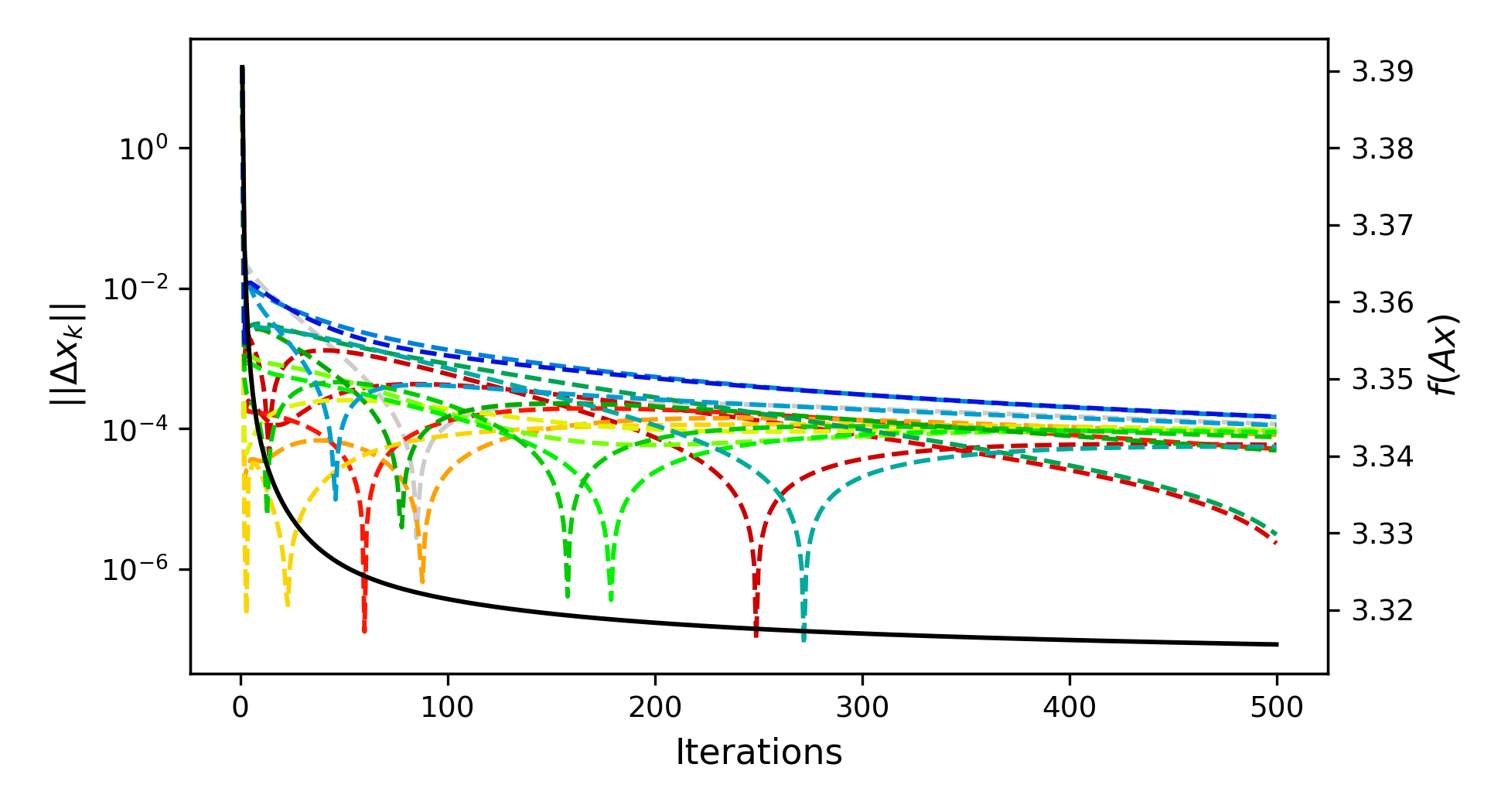}
    \caption[Stopping criteria for the inverter dataset]{Local updates (dashed) vs loss function (solid) for the inverter voltage dataset (top: original data, bottom: preprocessed data)}
    \label{fig:inv-stop}
\end{figure}

To observe the effect of data preprocessing on the proposed dynamic stopping criterion, the values $f(\Aa \xx)$ and $|\Delta \xk|$ are recorded for the first $500$ iterations and presented in Figure \ref{fig:inv-stop}.
With the original dataset, $f(\Aa \xx)$ exhibits a very slow linear decay:  $f(\Aa \xx) \sim \varepsilon n$, with $\varepsilon \ll 1$.
Preprocessing the data results in at least a power-law decay in $f(\Aa \xx)$.
Regarding using the values $|\Delta \xk|$ as a stopping criterion, preprocessing the data appears to be necessary to avoid stagnation in the decrease of values $|\Delta \xk|$, as was already seen for the decrease of
values $\Oa$.
At this point, it is worth noting that the data matrix of the inverters voltage dataset was found to be of low numerical rank: a Lasso regression leads to more zero regression coefficients than non-zero coefficients.
Therefore, the stagnation behavior can be explained by the \cola{} algorithm oscillating between two or more non-unique solutions of the least-square problem.
While the non-uniqueness might be removed or improved by further refinement of the elastic net regularization parameters $\lambda$ and $\eta$, the use of preprocessing avoids the need of parameter tuning and allows for the choice of $g=0$.

\begin{figure}[!htbp]
    \centering
    \includegraphics[width=0.8\linewidth]{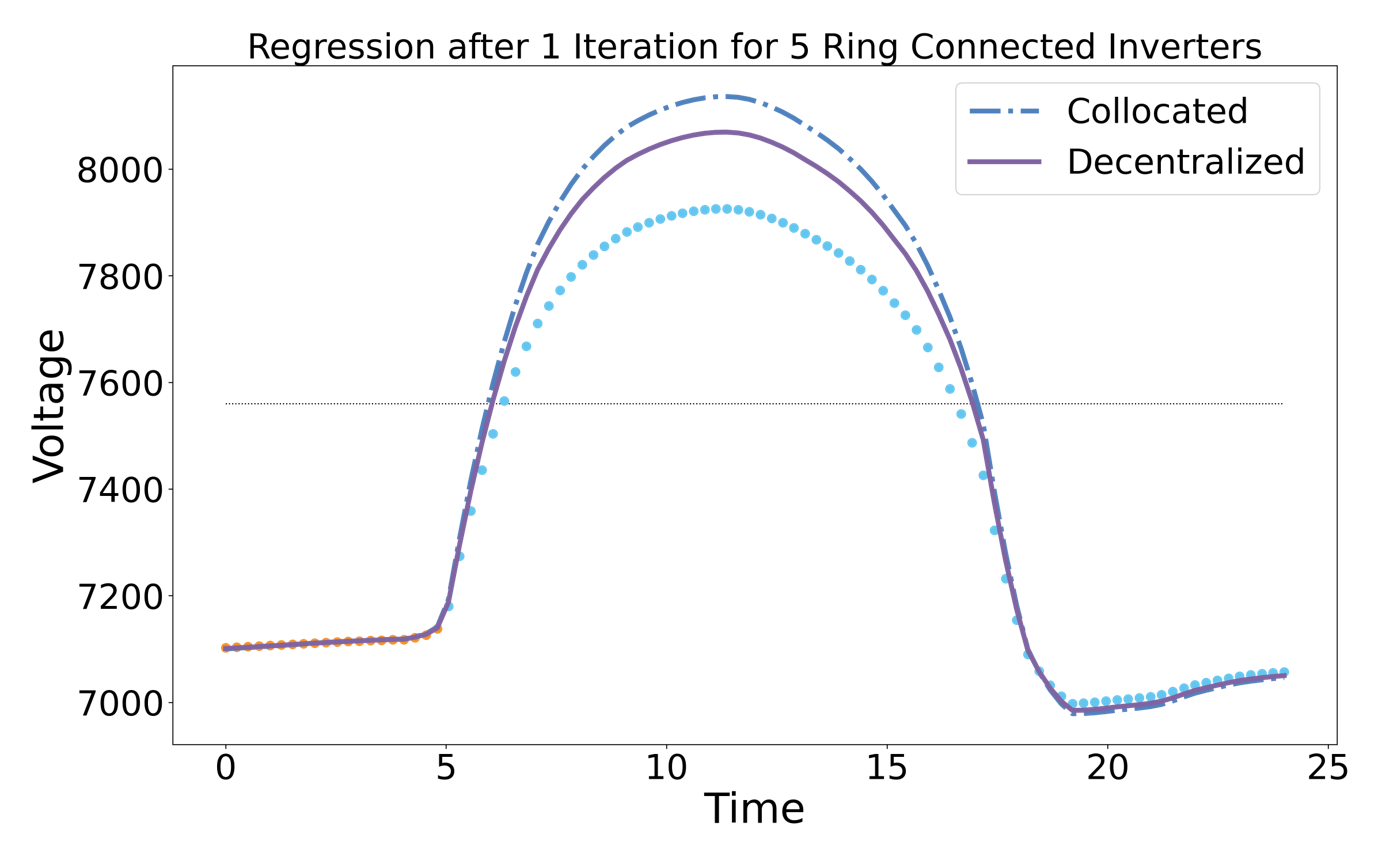}
    \caption[Overvoltage Scenario Prediction]{Comparison of Lasso regression with \cola{} (decentralized) prediciton and \texttt{sklearn} (collocated) prediction (orange dots: training data, blue dots: test data).}
    \label{fig:ov-scenario}
\end{figure}

Whereas the regression in Figure \ref{fig:inv-regression} is trained on a random
sample of the data, a more practical setting might be to train the regression on
a contiguous section of the data.  Consider, for example, a regression trained
on the first $5$ hours of data that is then used to predict voltages for the
latter $19$ hours.  Furthermore, an overvoltage is created by substantially
increasing the PV in the $80\%$ distribution network.  Two linear least-squares
models are trained and tested as to whether the admittedly exaggerated
overvoltage scenario is detected.  The first model, referred to as the
\textit{decentralized model}, is trained using a single communication round of
the \cola{} algorithm with a Lasso approach ($\eta = 0$).  The second model,
referred to as the \textit{collocated model}, is generated using a centralized
Lasso regression from the SciKit-Learn Python library \cite{scikit-learn}.
Figure \ref{fig:ov-scenario} shows that \cola{} outperforms the collocated
regression at accurately predicting both when the overvoltage scenario begins
and ends.

\section{Conclusion \label{sec:conclusion}}
For a particular distribution network augmented with $80\%$ PV penetration, a linear least-square regression is able to accurately predict $364$ days of inverter voltage magnitude data after only $1$ day of training.
The \cola{} algorithm is shown capable of producing such a linear least-square regression in a decentralized and asynchronous fashion.
Furthermore, a suitable regression for an overvoltage scenario is obtained by \cola{} with substantial reduction in communication cost compared to aggregating all the data in one location.
Finally, it was found that augmenting the \cola{} algorithm with a data preprocessing step is crucial to obtain sustained decay in both the objective function and the local update magnitudes.
An investigation into the extent that the results of this work apply to other distribution networks and a theoretical analysis how the preprocessing step affects the \cola{} algorithm are subjects of future work.

\section*{Disclaimer}
This document was prepared as an account of work sponsored by an agency of the
United States government.  Neither the United States government nor Lawrence
Livermore National Security, LLC, nor any of their employees makes any warranty,
expressed or implied, or assumes any legal liability or responsibility for the
accuracy, completeness, or usefulness of any information, apparatus, product, or
process disclosed, or represents that its use would not infringe privately owned
rights.  Reference herein to any specific commercial product, process, or
service by trade name, trademark, manufacturer, or otherwise does not
necessarily constitute or imply its endorsement, recommendation, or favoring by
the United States government or Lawrence Livermore National Security, LLC.  The
views and opinions of authors expressed herein do not necessarily state or
reflect those of the United States government or Lawrence Livermore National
Security, LLC, and shall not be used for advertising or product endorsement
purposes.

\bibliographystyle{IEEEtran}
\bibliography{refs}

\end{document}